\documentclass{l4dc2026}

\usepackage{algorithm}
\usepackage{algorithmic}
\usepackage{bm}

\usepackage{forest}
\DeclareMathOperator{\vecc}{vec}	    

\usepackage{mathtools}

\DeclarePairedDelimiter\norm{\lVert}{\rVert}

\newtheorem{assumption}{Assumption}
\newtheorem*{proofnew}{Proof}

\title[Subgradient Method for System Identification]{Subgradient Method for System Identification \\with Non-Smooth Objectives}
\usepackage{times}

\coltauthor{\Name{Baturalp Yalcin} \Email{byalcin@berkeley.edu}\\
 \Name{Jihun Kim} \Email{jihun.kim@berkeley.edu}\\
 \Name{Javad Lavaei} \Email{lavaei@berkeley.edu}\\
 \addr Department of Industrial Engineering and Operations Research,\\ University of California, Berkeley}

\begin{document}

\maketitle

\begin{abstract}%
This paper investigates a subgradient-based algorithm to solve the system identification problem for linear time-invariant systems with non-smooth objectives. This is essential for robust system identification in safety-critical applications. While existing work provides theoretical exact recovery guarantees using optimization solvers, the design of fast learning algorithms with convergence guarantees for practical use remains unexplored. We analyze the subgradient method in this setting, where the optimization problems to be solved evolve over time as new measurements are collected, and we establish linear convergence to the ground-truth system for both the best and Polyak step sizes after a burn-in period. We further characterize sublinear convergence of the iterates under constant and diminishing step sizes, which require only minimal information and thus offer broad applicability.  Finally, we compare the time complexity of standard solvers with the subgradient algorithm and support our findings with experimental results. This is the first work to analyze subgradient algorithms for system identification with non-smooth objectives.%
\end{abstract}

\section{Introduction}

Dynamical systems form the foundation of modern problems, such as reinforcement learning, auto-regressive models, and control systems. These systems evolve based on their current state, with the system dynamics equations determining the next state. Due to the complexity of real-world systems, these dynamics are often unknown or difficult to model precisely. Learning the unknown dynamics of a dynamical system is known as the \textit{system identification} problem.

The system identification problem has been extensively studied in the literature \citep{chen2012identification, lennart1999system}. The least-squares estimator (LSE) is the most widely studied approach for this problem as it provides a closed-form solution. Early research focused on the asymptotic (infinite-time) convergence properties of LSE \citep{ljung1978convergence}. Recently, finite-time (non-asymptotic) analyses gained popularity, leveraging high-dimensional statistical techniques for both linear \citep{simchowitz2018learning, faradonbeh2018finite, tsiamis2019finite, dean2018sample} and nonlinear dynamical systems \citep{foster2020learning, sattar2022non, ziemann2022single}. See \cite{tsiamis2023} for a survey. Traditionally, both asymptotic and non-asymptotic analyses of LSE focused on systems with independent and identically distributed (sub)-Gaussian disturbances. However, this assumption is often unrealistic for safety-critical systems such as autonomous vehicles, unmanned aerial vehicles, and power grids. In these settings, disturbances may be designed by an adversarial agent, leading to dependencies across time and arbitrary variations.

Since LSE is susceptible to outliers and temporal dependencies, recent research focused on robust non-smooth estimators for safety-critical systems. \citet{feng2021learning} and \citet{feng2023learning} analyzed the properties of these estimators through the lens of the Null Space Property. Later works modeled disturbances as zero with probability $1-p$ and nonzero with probability $p$.  For any $p$ between $0$ and $1$ with a directional symmetry, \citet{yalcin2024exact} and \citet{zhang2024exactrecoveryguaranteesparameterized} proved that non-smooth estimators can exactly recover system dynamics for linear and parameterized nonlinear systems, respectively. Moreover, \citet{kim2024prevailing} showed that asymmetric disturbances with nonzero probability $p$ could be transformed into directional-symmetric disturbances with probability $2p$, enabling exact recovery for any disturbance structure where $p < 1/2$. Building on these works, \cite{kim2025system} and \cite{kim2025sharp} respectively demonstrated the recovery of input-output mapping of partially observed systems and parameterized nonlinear systems. 
These non-asymptotic analyses leverage techniques from modern high-dimensional probability and statistics \citep{vershynin_2018, wainwright_2019}. Despite their strong learning properties, these approaches rely on convex optimization solvers, which become computationally prohibitive as the system state dimension and time horizon grow. Safety-critical systems require rapid estimation to ensure real-time control without delays. Thus, it is imperative to design fast algorithms that efficiently solve non-smooth estimators within a reasonable time frame. 

Due to the non-smoothness of the objective, gradient- and Hessian-based first-order and second-order algorithms are not viable in this context. Instead, we leverage the \textit{subgradient method}, a well-established technique for non-smooth convex optimization. Unlike gradient-based methods, subgradient methods do not guarantee descent at each iteration. Instead, existing results show that the minimum sub-optimality gap asymptotically converges to zero over time with both constant and diminishing step sizes. For a comprehensive overview of subgradient methods, see \cite{polyak1987introduction, bertsekas1997nonlinear, nesterov2009primal, 
shor2012minimization, nesterov2013introductory, kim2025revisit}. Closed-form solutions and gradient-based methods are widely used in system identification problems, as the literature has primarily focused on LSE or other smooth estimators for linear systems \citep{keesman2011system}. Nonlinear system identification is inherently more complex than its linear counterpart \citep{schoukens2019nonlinear}, yet existing estimation techniques still rely on smooth penalty functions. In contrast, subgradient methods have not been explored in the system identification context. Thus, this paper is the first to analyze the subgradient algorithm for the system identification problem. It is important to note that we do not apply the subgradient method to a single optimization problem but instead to a sequence of time-varying problems, since the optimization problem to be solved evolves as new measurements are sequentially taken. 
Within this framework, \textit{\textbf{our contribution is twofold}}: we show that, with high probability,  the algorithm iterates (1) \textit{converge linearly} to the ground-truth system under the best and Polyak step sizes after a burn-in period, and (2) \textit{converge sublinearly} to the ground-truth system under constant and diminishing step sizes after a burn-in period.

We note that, after the first appearance of this paper, \cite{kim2025bridge} extended our work to further reduce computational complexity using a stochastic subgradient method, where the expected estimation error is provably guaranteed to converge to the ground-truth system.

\textbf{Notation:} For a matrix $Z$, $\| Z \|_{F}$ denotes its Frobenius norm. For two matrices $Z_1$ and $Z_2$, we use $\langle Z_1, Z_2 \rangle$ to denote the inner product. 
For a vector $z$, $\| z \|_2$  denotes its $\ell_2$-norm. 
Given two functions $f$ and $g$, the notation $f(x) = \Theta(g(x))$ means that there exist universal positive constants $c_1$ and $c_2$ such that $c_1 g(x) \leq f(x) \leq c_2 g(x) $. Similarly, $f(x) = \mathcal{O}(g(x))$ implies that there exists $c_3>0$ such that $f(x) \le c_3 g(x) $, and $f(x) = \Omega(g(x))$ implies that there exists  $c_4>0$ such that $f(x)\geq c_4 g(x) $. $I_n$ stands for the $n \times n$ identity matrix.

\section{Non-smooth Estimator and Disturbance Model}

We consider a linear time-invariant (LTI) dynamical system of order $n$ with the update equation 
\begin{align}\label{sysdyn}
    x_{t+1} = \bar A x_t + \bar d_t, \quad t= 0, 1, \dots, T-1,
\end{align}
 where $\bar A \in \mathbb{R}^{n \times n}$ is the unknown system matrix and $\bar d_t \in \mathbb{R}^n$ are unknown system disturbances. Our goal is to estimate $\bar A$ using the state samples $\{ x_t \}_{t=0}^{T}$ obtained from a single system initialization. Note we only study the autonomous case due to space restrictions and the generalization to the case with inputs is straightforward. Traditionally, the following LSE has been widely studied in the literature:
\begin{equation}
    \begin{aligned}
        \min_{A \in \mathbb{R}^{n \times n}}\;\,  \sum_{t=0}^{T-1} \norm{x_{t+1} - A x_t}_2^2.           
    \end{aligned} \label{eq:lse} \tag{LSE}
\end{equation}
However, the LSE is vulnerable to outliers and adversarial disturbances.  In safety-critical system identification, a new popular alternative is the following non-smooth estimator:
\begin{equation}
    \begin{aligned}
        \min_{A \in \mathbb{R}^{n \times n}}\;\,  \sum_{t=0}^{T-1} \norm{x_{t+1} - A x_t}_2.           
    \end{aligned} \label{eq:hard-lasso2} \tag{NSE}
\end{equation}

The primary objective of this non-smooth estimator is to achieve exact recovery in finite time, which is an essential feature for safety-critical applications such as power-grid monitoring, unmanned aerial vehicles, and autonomous vehicles. In these scenarios, the system should be learned from a single trajectory (since multiple initialization is often infeasible) and, moreover, it can be exposed to adversarial agents. Even under i.i.d. disturbance vectors, if $\bar d_t$ is non-zero at every time step $t$, the exact recovery of the $\bar A$ is not attainable in finite time. Motivated by \cite{zhang2024exactrecoveryguaranteesparameterized}, we, therefore, assume a disturbance model in which the disturbance vector $\bar d_t$ is zero at certain time instances. 
\begin{assumption}[Probabilistic sparsity model]\label{def:probabilistic}
For each time instance $t$, the disturbance vector $\bar{d}_t$ is zero with probability $1-p$, where $p \in (0,1)$. 
Furthermore, the time instances at which the disturbances are nonzero are independent of each other (while the nonzero disturbance values could be correlated).
\end{assumption}

 It prevails the system identification literature that persitent excitation is crucial to successfully learn the system \citep{NarendraAnnaswamy1987}. As a result, lots of literatures  in  \eqref{eq:lse} also adopt the concept of persistent excitation and so do us. For example, \cite{zhang2024exactrecoveryguaranteesparameterized} showed in their Theorem 3 that it is infeasible to learn the full system dynamics
if the disturbance vectors take arbitrary values and are chosen from a low-dimensional subspace.
 To this end, we define the filtration $\mathcal F_t = \bm{\sigma}\{x_0, \dots, x_t\}$ and present a persistent excitation condition, which ensures exploration of the entire state space.
 
\begin{assumption}[Persistent Excitation]\label{persistent}
   Conditional on the past information $\mathcal{F}_t = \bm{\sigma}\{x_0, \dots, x_t\}$ and the event that $\bar d_t \neq 0$, the disturbance vector satisfies
   \begin{align}
       \mathbb{E}[(\bar A x_t + \bar d_t)(\bar A x_t + \bar d_t)^T~|~\mathcal{F}_t, ~\bar d_t \neq 0] \succeq\lambda^2 I_n,\quad \forall t=0,\dots, T-1,
   \end{align}
    where $\lambda >0$ is a constant.    
\end{assumption}

We now impose a sufficient condition for the system stability of the LTI system below to guarantee learning under possibly large correlated disturbances.
\begin{assumption}
    \label{asp:stable}
The ground-truth $\bar{A}$ satisfies $\rho := \left\|\bar{A}\right\|_2 < 1$ (note that this assumption can easily be extended to the stability of $\bar A$—spectral radius being less than 1—due to Gelfand's formula). 
\end{assumption}

We also assume below that the disturbance vectors $\bar d_t$ are sub-Gaussian, which limits the tail behavior of $\bar d_t$ (note that bounded disturbances are automatically sub-Gaussian). However, they are temporally correlated since each disturbance $\bar d_t$ may depend on previous disturbances. This correlation better reflects real-world attacks, where an adversarial agent can inject disturbances that depend on prior injections.

\begin{assumption}[Direction-restricted sub-Gaussian disturbances]\label{asp:subgaussian}
Conditional on the filtration $\mathcal{F}_t$ and the event that $\bar{d}_t\neq 0$, 
the attack vector $\bar{d}_t$ is defined by the product $\ell_t  \hat{d}_t$, where 
\begin{enumerate}
    \item  $\hat d_t$ is a zero-mean unit vector whenever $\bar{d}_t \neq 0$, for any $\mathcal{F}_t$;
    \item $\ell_t$ is a finite, positive variable that the adversary can arbitrarily select;
    \item $\bar d_t = \ell_t \hat d_t$ is sub-Gaussian with parameter $\sigma$, for any $\mathcal{F}_t$. 
\end{enumerate}
\end{assumption}

\begin{remark}
In Appendix \ref{app:signrestriction}, we show that Assumptions \ref{def:probabilistic}, \ref{persistent}, \ref{asp:stable} and \ref{asp:subgaussian} suffice to derive the same finite sample complexity bounds for exact recovery in \cite{zhang2024exactrecoveryguaranteesparameterized}.  
Note that $\hat d_t$ and $\ell_t$ implies a direction and a length of the disturbance, respectively. 
    At time $t$, the adversary can select any positive $\ell_t$ that is not necessarily independent of $\hat d_t$, and even select the maximally adversarial $\ell_t$ based on the previous information $\mathcal{F}_t$. However, one may notice that the restriction on the direction $\hat d_t$  limits the behavior of disturbances. In fact, Theorem 3 in \cite{kim2024prevailing} showed that when the probabilistic sparsity parameter $p\in (0,1/2)$, this restrictive assumption can be removed and an arbitrary distribution $\bar d_t$ is approximately within the class of disturbances defined in Assumption \ref{asp:subgaussian}. This incorporates realistic scenarios where the attack happens infrequently but the magnitude of attacks can be arbitrary. 
\end{remark}

In recent years, non-smooth estimators for system identification gained traction \citep{feng2021learning, feng2023learning, yalcin2024exact, zhang2024exactrecoveryguaranteesparameterized, kim2024prevailing, kim2025system, kim2025sharp}. These studies analyzed non-smooth estimators in a non-asymptotic framework under various disturbance models. However, the aforementioned works do not propose numerical algorithms and rely on existing convex optimization solvers, which suffer from two drawbacks. The first issue is that in practice \eqref{eq:hard-lasso2} should be repeatedly solved for $T=1,2,...$ as new measurements are collected until a satisfactory model is obtained. This means that the methods in the aforementioned works require solving a convex optimization problem at each time step. In addition, the problem \eqref{eq:hard-lasso2} at time $T$ is translated into a second-order conic program with $T + n^2$ scalar variables and $T$ conic inequalities. Hence, the solver performance degrades significantly as the order of the system and the time horizon grow. Since safety-critical applications demand fast and efficient identification, convex optimization solvers are inefficient for real-time computation. To address these limitations, we investigate the subgradient method for learning system parameters.

\section{Subgradient Method}

We let $f_T(A)$ denote the objective function of \eqref{eq:hard-lasso2}; \textit{i.e.}, $f_T(A)=\sum_{t=0}^{T-1} \norm{x_{t+1} - A x_t}_2$, and $\partial f_T(A)$ denote its subdifferential at a point $A$ and a time period $T$. We denote a subgradient in the subdifferential by $G_{A, T} \in \partial f_T(A)$. We first define the subdifferential of the $\ell_2$ norm.
\begin{definition} The subdifferential of the $\ell_2$ norm is defined as
$\partial \| x \|_2 = \begin{cases}
        \frac{x}{\|x\|_2}, & \text{if } x \not = 0, \\
        \{ e :  \|e\|_2 \le 1 \}, & \text{if } x=0.
    \end{cases}$
\label{def:subgradient}
\end{definition}
The subdifferential $\partial \|x\|_2$ is a unique point when $x$ is nonzero, but it consists of all points within the unit ball when $x=0$. Using the previous definition, we obtain the subdifferential of the objective function, $\partial f_T(A)$.
\begin{lemma}\label{subd_of_f_T} The subdifferential of $f_T(\cdot)$ at a point $A$ is
     $   \partial f_T(A) = \Bigr\{  - \sum_{t = 0}^{T-1} \partial \| x_{t+1} - A x_t \|_2   \cdot x_t^T   \Bigr\},$
and the subdifferential of the $f_T(\cdot)$ at the point $\bar A$ is
    \begin{align*}
        \partial f_T(\bar A) = \Biggr\{ - \sum_{t \in \mathcal{K}} \hat d_t  x_t^T  - \sum_{t \in \mathcal{K}^c} e_t  x_t^T : \|e_t\|_2 \le 1, \forall t \in \mathcal{K}^c\Biggr\},
    \end{align*}
where $\mathcal{K} := \{t\in \{0, \dots, T-1\} ~|~ \bar{d}_t \neq 0\}$ is the set of time indices of nonzero disturbances,  
$\mathcal{K}^c :=\{0,\dots, T-1\}\setminus \mathcal{K} $ is the complement of $\mathcal{K}$, 
and $\hat{d}_t := \bar{d}_t /  \|\bar{d}_t\|_2,\forall t\in\mathcal{K}$, are the normalized disturbance directions. \label{lem:min}
\end{lemma}
The proof of Lemma \ref{lem:min} simply follows from a property of the Minkowski sum over convex sets. 

%
To find the matrix $\bar A$, we propose Algorithm \ref{alg} that relies on the subgradient method. In Step (i), the subdifferential of the objective function at the current estimate is calculated. Since this subdifferential possibly contains multiple values, we choose a subgradient from the subdifferential in Step (ii). After that, a step size is selected in Step (iii). In Step (iv), a subgradient update is performed. Lastly, we update the objective function with the new data point in Step (v). The main feature of this algorithm is that we do not apply the subgradient algorithm to solve the current optimization problem \eqref{eq:hard-lasso2} and instead we take one iteration to update our estimate of the matrix value at time $T$. Then, after a new measurement is taken, the optimization problem at time $T+1$ changes since an additional term is added to its objective function. Thus, as we update our estimates, the optimization problems also change over time. We analyze the convergence of Algorithm \ref{alg} under different step size selection rules and show that although the objective function changes at each iteration, our algorithm converges to the $\bar A$.

\begin{algorithm}[ht]
    \caption{Subgradient Method for \eqref{eq:hard-lasso2}}\label{alg}
    \begin{algorithmic}
        \STATE Let $f_1(A) = \|x_1-Ax_0\|_2$.  Initialize $\hat A^{(1)}$ randomly.
        \FOR{$T = 1, 2,\dots$}
            \STATE (i) Find $\partial f_T(\hat A ^{(T)})= - \sum_{t=0}^{T-1} \partial \norm{(\bar A -\hat A^{(T)}) x_t + \bar d_t}_2 \cdot x_t^T$.
            \STATE (ii) Choose $G_{\hat A^{(T)}, T}$ from $\partial f_T(\hat A ^{(T)})$.
            \STATE (iii) Choose $\beta^{(T)}$ using an appropriate step size rule.
            \STATE (iv) Update $\hat A^{(T+1)} = \hat A^{(T)} -  \beta^{(T)}G_{ \hat A^{(T)}, T}$.
            \STATE (v) Update  $f_{T+1}(A) = f_T(A) + \| x_{T+1} - Ax_{T}\|_2.$
        \ENDFOR
        \STATE Return $\hat A^{(2)},\hat A^{(3)},\hat A^{(4)}, \dots $
    \end{algorithmic}
\end{algorithm}

\subsection{Subgradient Method with Best Step Size}

We analyze the subgradient method using the best step size, which is defined as the step size that maximizes improvement at each iteration. However, computing this step size requires knowledge of the ground-truth matrix $\bar A$, making it impractical in real-world applications. Nonetheless, studying this case offers valuable insights into the behavior of subgradient updates. Let the step size in Algorithm \ref{alg} be chosen as 
\[ \beta^{(T)} =\frac{ \langle G_{\hat A^{(T)}, T}, \hat A^{(T)} - \bar A \rangle}{\|G_{\hat A^{(T)}, T}\|_F^2}.\]
It is shown in Appendix \ref{app:proof} that this is the best step size possible. According to Theorem \ref{thm:angle-cosine}, the distance to the ground-truth matrix decreases by a factor of $\sqrt{1 - \cos^2({\hat \theta_T})} $ at each update with the best step size. As long as this angle is not $\pm 90^\circ$, Algorithm \ref{alg} makes progress toward $\bar A$.

\begin{theorem} \label{thm:angle-cosine}
    Under the best step size for Algorithm \ref{alg}, it holds that
    \[ \| \hat A^{(T+1)} - \bar A \|_F \le \sqrt{1 - \cos^2({\hat \theta_T})} \| \hat A^{(T)}- \bar A \|_F,\quad \forall T\geq 1,\]
    where $\hat \theta_T$ denotes the angle between $\hat A^{(T)} - \bar A$ and $G_{\hat A^{(T)}, T}$.
\end{theorem}
 If $\hat \theta_T$ is obtuse, the step size will be negative, which must be avoided. To this end, we establish a positive lower bound on $\cos(\hat \theta_T)$ in Theorem \ref{thm:cosine-lower-bound} and show that this angle remains acute. To guarantee that, we define the burn-in period as the number of samples required for a theoretical exact recovery. With high probability, $\bar A$ turns out to be the unique optimal solution of $\min_A f_T(A)$ for all time steps after the burn-in period.

\begin{definition} \label{def:burn-in}
    Consider an arbitrary $\delta \in (0,1)$ and constants $p, \lambda, \rho, \sigma$ each from Assumptions \ref{def:probabilistic}, \ref{persistent}, \ref{asp:stable} and \ref{asp:subgaussian}. Let $\kappa:= \sigma/\lambda$. Let $T^{(burn)}$ be defined as 
    \begin{align*}\label{eqn:lipschitz-upper}
         T^{(burn)}& :=  \Theta\biggr(\max\left\{  \frac{\kappa^{10}}{(1-p)^2(1-\rho)^3 }, \frac{\kappa^4}{p(1-p)}\right\}\cdot \left[ n^2 \log\left(\frac{\kappa}{p(1-p)(1-\rho)}\right)+\log\left(\frac{1}{\delta}\right)\right]\biggr).
    \end{align*}%
Then,  with probability at least $1-\delta$, the ground-truth matrix $\bar A$ is the unique optimal solution to \eqref{eq:hard-lasso2} for all $T \ge T^{(burn)}$  (see Appendix \ref{app:signrestriction} for details). 
\end{definition}
%
We analyze the algorithm's behavior after the burn-in period because it happens that the iterates of Algorithm \ref{alg} are chaotic when $T$ is not large enough or equivalently when $\bar A$ is not the unique solution of \eqref{eq:hard-lasso2}. The following theorem establishes that if the number of samples exceeds $T^{(burn)}$, the cosine of the angle can be bounded below by a positive term. 

 \begin{theorem} \label{thm:cosine-lower-bound}
    Define $\hat \theta_T$ to be the angle between $ \hat A^{(T)} - \bar A$ and $G_{\hat A^{(T)}, T}$.   Suppose that  Assumptions 
   \ref{def:probabilistic}, \ref{persistent}, \ref{asp:stable} and \ref{asp:subgaussian} hold. Then, with probability at least $1-\delta$, the term $\cos({\hat \theta_T})$ can be bounded below by 
    \[ \cos({\hat \theta_T}) =\Omega\left(\frac{  p(1-p) (1-\rho)  }{ \kappa^5  }  \right)> 0.\]
for all $T \geq T^{(burn)}$ such that  $\hat A^{(T)} \neq \bar A$.
\end{theorem}
The analysis of Theorem \ref{thm:cosine-lower-bound} exploits the sharpness of 
 $f_T(\cdot)$ with respect to $\bar A$ when $T\ge T^{(burn)}$; \textit{i.e.}, $\exists \bar c>0$ such that $f_T(\hat A^{(T)})-f_T(\bar A) \geq \bar c \|\hat A^{(T)}- \bar A\|_F$
 (see Appendix \ref{pf:cosine} for details).  Theorem \ref{thm:cosine-lower-bound} provides a probabilistic lower bound on the cosine of the angle of interest when the sample complexity $T$ is sufficiently large. 
Together with Theorem \ref{thm:angle-cosine}, it can be readily shown that the algorithm iterates converge toward the ground-truth $\bar A$ after a sufficiently large time $T^{(burn)}$. 

%
\begin{corollary} \label{cor:linear-convergence-best}
  Suppose that  Assumptions 
   \ref{def:probabilistic}, \ref{persistent}, \ref{asp:stable} and \ref{asp:subgaussian} hold.   Let $D = \| \hat A^{(T^{(burn)})} - \bar A \|_F$ denote the distance of the solution of the last iterate to the ground-truth $\bar A$ at the end of the burn-in time. Define $\gamma$ to be 
    \[ \gamma := \sqrt{1-  \Omega\left(\frac{p^2(1-p)^2(1-\rho)^2}{\kappa^{10}}\right) }.\]
    Then, given $\delta\in(0,1)$ and an arbitrary $\epsilon > 0$, with probability at least $1-\delta$,     
    Algorithm \ref{alg} under the best step size achieves $\| \hat A^{(T)} - \bar A \|_F \le \epsilon$  when the number of iterations satisfies $T \ge  T^{(burn)} + T^{(conv)}(\epsilon)$, where
    \begin{align*}
        T^{(conv)}(\epsilon) := \frac{\log\left( D/\epsilon\right)}{-\log(\gamma)}.
    \end{align*}
\end{corollary}
The proof simply follows from identifying $T$ such that $D\cdot \gamma^T \le \epsilon$. 
Corollary \ref{cor:linear-convergence-best} implies that to achieve an $\epsilon$-level estimation error, the total number of required samples is the sum of those needed for the burn-in phase $T^{(burn)}$ and the linear convergence phase $T^{(conv)}(\epsilon)$. With high probability, the burn-in period scales as $\max\{(1-p)^{-2}, p^{-1}(1-p)^{-1} \}$ with respect to $p$ and as $n^2$ with respect to the system dimension. After the burn-in phase, the algorithm converges to the ground-truth system at a linear rate because $\cos(\hat \theta_T)$ is bounded below by a positive factor, and thus the linear convergence rate is bounded above by $\gamma<1$. This ensures that the sample complexity scales logarithmically as $\log(D/\epsilon)$ in the convergence phase. 
In the next subsections, we now generalize the results to more practical step size rules. 

\subsection{Subgradient Method with Polyak Step Size}

In this subsection, we analyze the Polyak step size. Unlike the best step size selection, the Polyak step size does not require knowledge of $\bar A$. However, it requires information about the optimal objective value at $\bar A$ at time $T$; \textit{i.e.}, $f_T(\bar A)$. Under this step size rule, we can establish the result similar to Corollary \ref{cor:linear-convergence-best}.
\begin{theorem} \label{thm:polyak-result}
   Suppose that  Assumptions 
   \ref{def:probabilistic}, \ref{persistent}, \ref{asp:stable} and \ref{asp:subgaussian} hold. Let $\delta\in(0,1)$ and $\epsilon > 0$ be arbitrary. Consider Algorithm \ref{alg} under the Polyak step size defined as 
    \begin{align*}
        \beta^{(T)} = \frac{f_T(\hat A^{(T)}) - f_T(\bar A) }{\| G_{\hat A^{(T)}, T} \|_F^2}.
    \end{align*}
With probability at least $1-\delta$,     
    Algorithm \ref{alg} under the Polyak step size achieves $\| \hat A^{(T)} - \bar A \|_F \le \epsilon$  when $T \ge  T^{(burn)} + T^{(conv)}(\epsilon)$, where $T^{(conv)}(\epsilon)$ is defined in Corollary \ref{cor:linear-convergence-best}.
\end{theorem}
The proof of Theorem \ref{thm:polyak-result} can be found in Appendix \ref{pf:polyak}. Although the Polyak step size relies on less information about the ground-truth system dynamics than the best step size (\textit{i.e.}, it needs to know an optimal scalar rather than an optimal matrix), it requires the same number of iterations as the best step size selection. These step sizes are implementable in practice (subject to some approximations) if some prior knowledge about the system dynamics is available.

\subsection{Subgradient Method with Constant and Diminishing Step Sizes}

The step sizes discussed in previous subsections require knowledge of the ground-truth matrix $\bar A$ or the optimal objective value $f_T(\bar A)$. While theoretically useful, these step sizes are impractical in real-world settings. We therefore explore the subgradient method with constant and diminishing step sizes, which are broadly applicable. Unlike the best and Polyak step sizes, constant and diminishing step sizes do not guarantee descent at every iteration. Instead, we analyze the \textit{minimum solution gap},  $\min_{T\in[T^{(burn)},\bar T]} \|\hat A^{(T)}-\bar A\|_F$, and show that it converges to zero when $\bar T$ is sufficiently large.

Now, we demonstrate that under a constant step size, the minimum solution gap of the subgradient algorithm polynomially converges at a rate of $\bar T^{-1/2}$, thereby showing sublinear convergence. Specifically, Theorem \ref{thm:constant-step-size} states that when $\bar T$ is significantly larger than $T^{(burn)}$,  the minimum solution gap can be simplified as $\mathcal{O}\bigr(\frac{D\kappa^5}{p(1-p)(1-\rho) \bar T^{1/2}}\bigr) $. 

\begin{theorem} \label{thm:constant-step-size}
Suppose that  Assumptions 
   \ref{def:probabilistic}, \ref{persistent}, \ref{asp:stable} and \ref{asp:subgaussian} hold.   Let $D = \| \hat A^{(T^{(burn)})} - \bar A \|_F$, and fix a sufficiently large $\bar T$.   Consider Algorithm \ref{alg} with the step size $\beta^{(T)}$ set to the constant value
    \begin{align*}
        \beta^{(T)} = \beta =  \Theta\left(\frac{D(1-\rho)}{\sigma  (\sum_{t=T^{(burn)}}^{\bar T} t^2)^{1/2}}\right),\quad \forall T\ge T^{(burn)}.
    \end{align*}
    Then, with probability at least $1- \delta$, the minimum solution gap at time $\bar T$ is bounded  as 
     \begin{align*}
        \min_{T\in [T^{(burn)},\bar T]} \|\hat{A}^{(T)} - \bar A\|_F  =\mathcal{O}\left(\frac{D\kappa^5}{p(1-p)(1-\rho)}  \cdot \frac{(\sum_{t=T^{(burn)}}^{\bar T} t^2)^{1/2}}{\sum_{t=T^{(burn)}}^{\bar T} t}\right).
    \end{align*}

\end{theorem}
 Furthermore, we establish a similar result for the diminishing step size, which is chosen to satisfy $\sum_T \beta^{(T)} \xrightarrow{} \infty$ and $\sum_T (\beta^{(T)})^2 < \infty$, ensuring sufficient exploration and stable convergence without excessive oscillations. Theorem \ref{thm: diminishing-step-size} establishes that, under a diminishing step size, the minimum solution gap polynomially converges at a rate of $\bar T^{-1/2}$, which also yields sublinear convergence.

\begin{theorem} \label{thm: diminishing-step-size}
  Under the same setting as Theorem \ref{thm:constant-step-size}, suppose that the step size $\beta^{(T)}$ takes the diminishing form $\beta^{(T)} = \beta/T$ for all $T\ge T^{(burn)}$, where $\beta$ is defined as
    \begin{align*}
        \beta := \Theta \left( \frac{D  (1-\rho)}{\sigma  (\bar T- T^{(burn)})^{1/2}} \right).
    \end{align*}
    Then, with probability at least $1- \delta$, the minimum solution gap at time $\bar T$ is bounded as 
    \begin{align*}
        \min_{T\in [T^{(burn)},\bar T]} \|\hat{A}^{(T)} - \bar A\|_F =  \mathcal{O}\left(\frac{D\kappa^5}{p(1-p)(1-\rho)}  \cdot \frac{1}{(\bar T-T^{(burn)})^{1/2}}\right).
    \end{align*}
\end{theorem}
The analysis leverages the sharpness of $f_T(\cdot)$ used in Theorem \ref{thm:cosine-lower-bound} (see Appendix \ref{app:constant-diminishing} for details).

\section{Time Complexity Analysis}

We compare the time complexities of the subgradient method and the exact solution of \eqref{eq:hard-lasso2} using CVX. For the subgradient algorithm,  each time period involves computing the subgradient, which has a per-iteration cost of $\mathcal{O}(Tn^2)$. Running the algorithm for $T$ iterations results in an overall complexity of $\mathcal{O}(T^2n^2)$. On the other hand, solving \eqref{eq:hard-lasso2} using CVX reduces the problem to a second-order conic program (SOCP) with $T + n^2$ variables with $T$ second-order cone constraints (minimization of the sum of norms can be written as a SOCP \citep{alizadeh2003second}). Defining $[T] = \{0, 1, \dots, T-1\}$, this problem is formulated as 
\begin{align} \label{eq: socp} 
    \min_{\substack{\vecc(A) \in \mathbb{R}^{n^2} \\ f_t \in \mathbb{R}, \; t \in [T]}}  \quad  \sum_{t=0}^T f_t \quad
     \text{s.t.}\quad   f_t \ge \| x_{t+1} - (x_t \otimes I_n) \vecc(A) \|_2, \; t \in [T], \tag{SOCP}
\end{align}
where $\otimes$ denotes the Kronecker product and $\vecc(A)$ is the vectorization of $A$ obtained by stacking its columns into a single vector. 
The CVX solver employs primal-dual interior point methods to solve \eqref{eq: socp}. The number of iterations required to reach an $\epsilon$-accurate solution is $\mathcal{O}(T^{1/2}\log(1/\epsilon))$ \citep{wright1997primal}. Since solving conic constraints has a worst-case complexity of $\mathcal{O}(n^6)$, the per time period complexity is $\mathcal{O}(T^{1/2}n^6)$, leading to an overall complexity of $\mathcal{O}(T^{3/2}n^6)$. We know that $T$ and $T^{(burn)}$ scale as $\Theta(n^2)$. Consequently, the per-time period complexity of the subgradient method and the exact SOCP solution are $\mathcal{O}(n^4)$ and $\mathcal{O}(n^7)$, respectively. This highlights that solving \eqref{eq: socp} exactly requires significantly more computational power than the subgradient approach.

\section{Numerical Experiments} \label{sec:experiments}

In this section, we present numerical experiments to validate the convergence results.
We generate a system with ground-truth matrix $\bar A\in \mathbb{R}^{n\times n}$ whose singular values are uniform in $(0,1)$. Given the probability $p$, disturbance vectors are set to zero with probability $1-p$. Then, we sample $ \ell_t \sim \mathcal{N}(0, \sigma_t^2)$, where $\sigma_t^2 := \min\{\|x_t\|_2^2, 1 / n\}$, and sample $\hat{d}_t$ from a uniform distribution over a unit sphere. 
When an attack occurs, the disturbance is set to $\bar{d}_t:= |\ell_t| \hat{d}_t$, which conforms to Assumption \ref{asp:subgaussian} while the disturbance values are still correlated over time. We generate $10$ independent trajectories of the subgradient algorithm using $10$ different systems. Then, we report the average solution gap and the average loss gap of Algorithm \ref{alg} at each period. The \textit{solution gap} and the \textit{loss gap} for Algorithm \ref{alg} at time $T$ are $\| \hat A^{(T)} - \bar A \|_F$ and $f_T(\hat A^{(T)}) - f_T(\bar A)$, respectively.

\begin{figure}[t]
    \centering
     \subfigure[\small Different Step Sizes for $n = 5$, and $p = 0.7$]{\label{fig:stepsizes}\includegraphics[height=83pt]{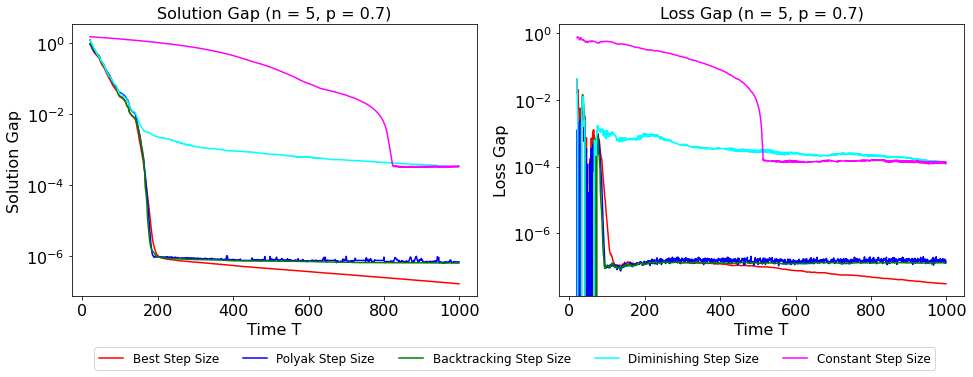}}\hfill
    \subfigure[\small Best Step Size, $p=0.7$,  $n\in \{50,75,100\}$]{\label{fig:longer}\includegraphics[height=83pt]{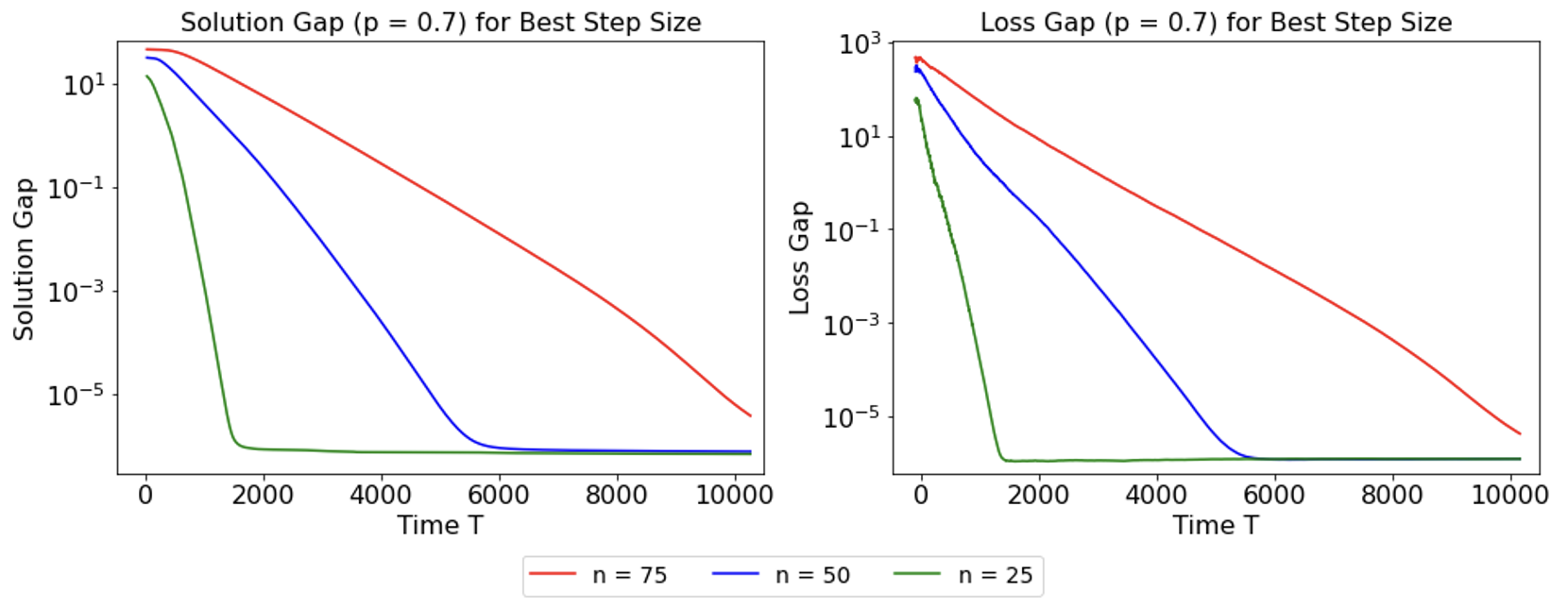}}
    \caption{\small Solution Gap and Loss Gap of Algorithm \ref{alg} with Different Step Sizes and  Dimensions $n$.}
    \label{fig:experiments}
\end{figure}

Algorithm \ref{alg} is experimented with different step size rules.  In addition to the best, Polyak, diminishing, and constant step sizes, we also introduce the backtracking (Armijo) step size, where the initial step size is repeatedly multiplied by a factor $\alpha \in (0,1)$ until we achieve a sufficient descent (see \cite{armijo1966minimization}). Figure \ref{fig:stepsizes} depicts the performance of Algorithm \ref{alg} under these different step sizes for $n=5$ and $p=0.7$. Best and Polyak step sizes achieve fast and linear convergence as proven theoretically.  Moreover, despite slower rates, the diminishing and constant step sizes converge to a neighborhood of the ground-truth system $A$. Interestingly, the backtracking step size performs similarly to the best and Polyak step sizes even without knowledge of $\bar A$ or $f_T(\bar A)$. Thus, the theoretical analysis of the backtracking step size would be a promising future direction.

\begin{figure}[t]
    \centering
     \subfigure[\small $p = 0.7, ~n \in \{5, 10, 15\}$]{\label{fig:dimensions}\includegraphics[height=87.5pt]{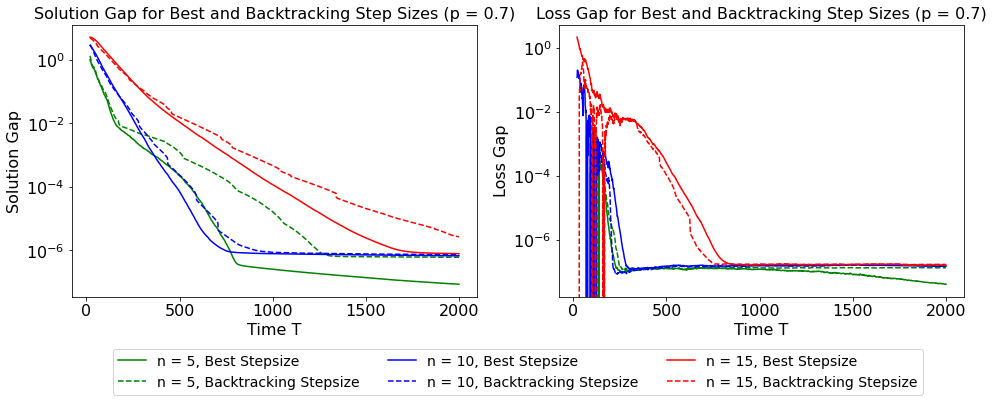}}\hfill
    \subfigure[\small  $n=5,  ~p \in \{0.5, 0.7, 0.8 \}$]{\label{fig:probabilities}\includegraphics[height=87.5pt]{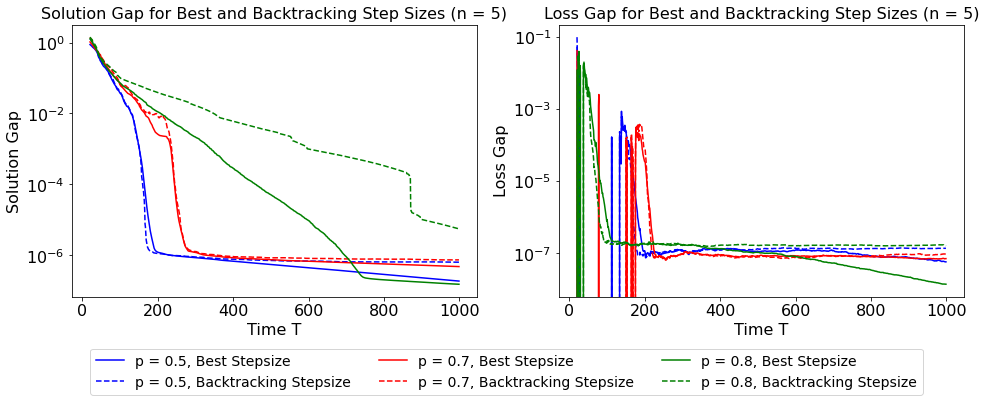}}
    \caption{\small Solution Gap and Loss Gap of Algorithm \ref{alg} with Best and Backtracking Step Sizes,  Different  Dimensions $n$ and Probabilistic Sparsity Parameters $p$.}
    \label{fig:experiments2}
\end{figure}

We test the performance of Algorithm \ref{alg} with different $n$ and $p$. We choose backtracking and best step sizes to test the impact of the system dimension, $n \in \{5, 10, 15 \}$, on the solution gap and the loss gap when $p=0.7$. Figure \ref{fig:dimensions} shows that as $n$ grows, the number of samples needed for convergence grows as $n^2$. This aligns with the theoretical results from earlier sections. Moreover, the backtracking step size performs similarly to the best step size. Likewise, Figure \ref{fig:probabilities} depicts the performance of the subgradient algorithm for  different probabilities of nonzero disturbance $p \in \{0.5, 0.7, 0.8 \} $ with $n = 5$. As $p$ increases, the convergence for the estimator needs more samples, as supported by the theoretical results.

%
Last but not least, we provide simulation results using the best step size for large systems of order $n \in \{ 25, 50, 75 \}$. We run these experiments for $T = 10000$. This is not viable with the standard CVX solver due to the high computational complexity of \eqref{eq: socp}. Nevertheless, the subgradient algorithm provides an accurate estimation for large systems within a reasonable time frame, as seen in Figure \ref{fig:longer}.

\section{Conclusion}

We propose and analyze a subgradient algorithm for non-smooth system identification of linear-time invariant systems with disturbances following a probabilistic sparsity model. While the non-asymptotic theory for this setting is well-established in the literature, existing methods rely on solvers that use primal-dual interior-point methods for second-order conic programs, leading to high computational costs. To address this, we propose a subgradient method that enables fast exact recovery. We demonstrate that the subgradient method with the best and Polyak step sizes converges linearly to the ground-truth system dynamics after a burn-in period. Moreover, we show that under both constant and diminishing step sizes, the minimum solution gap converges to zero at a polynomial rate, which is of greater practical relevance.
 These theoretical results are supported by numerical simulations.
Our experiments also highlight the practical effectiveness of the backtracking step size, suggesting it as a promising alternative. However, its theoretical properties remain an open question. Moreover, our current implementation computes the full subgradient at each iteration. To further improve time complexity, exploring the convergence behavior of partial subgradient approaches, such as mini-batch or stochastic gradient descent, within the system identification context would be an important theoretical direction.

\acks{This work was supported by the U. S. Army Research Laboratory and the U. S. Army Research Office under Grant W911NF2010219, Office of Naval Research under Grant N000142412673, AFOSR, NSF, and the UC Noyce Initiative.}

\bibliography{main}

\begin{thebibliography}{34}
\providecommand{\natexlab}[1]{#1}
\providecommand{\url}[1]{\texttt{#1}}
\expandafter\ifx\csname urlstyle\endcsname\relax
  \providecommand{\doi}[1]{doi: #1}\else
  \providecommand{\doi}{doi: \begingroup \urlstyle{rm}\Url}\fi

\bibitem[Alizadeh and Goldfarb(2003)]{alizadeh2003second}
Farid Alizadeh and Donald Goldfarb.
\newblock Second-order cone programming.
\newblock \emph{Mathematical Programming}, 95\penalty0 (1):\penalty0 3--51, 2003.

\bibitem[Armijo(1966)]{armijo1966minimization}
Larry Armijo.
\newblock Minimization of functions having lipschitz continuous first partial derivatives.
\newblock \emph{Pacific Journal of mathematics}, 16\penalty0 (1):\penalty0 1--3, 1966.

\bibitem[Bertsekas(1997)]{bertsekas1997nonlinear}
Dimitri~P Bertsekas.
\newblock Nonlinear programming.
\newblock \emph{Journal of the Operational Research Society}, 48\penalty0 (3):\penalty0 334--334, 1997.

\bibitem[Boyd and Vandenberghe(2004)]{boyd2004convex}
Stephen~P Boyd and Lieven Vandenberghe.
\newblock \emph{Convex optimization}.
\newblock Cambridge university press, 2004.

\bibitem[Chen and Guo(2012)]{chen2012identification}
Han-Fu Chen and Lei Guo.
\newblock \emph{Identification and stochastic adaptive control}.
\newblock Springer Science \& Business Media, 2012.

\bibitem[Dean et~al.(2020)Dean, Mania, Matni, Recht, and Tu]{dean2018sample}
Sarah Dean, Horia Mania, Nikolai Matni, Benjamin Recht, and Stephen Tu.
\newblock On the sample complexity of the linear quadratic regulator.
\newblock \emph{Foundations of Computational Mathematics}, 20\penalty0 (4):\penalty0 633--679, 2020.
\newblock \doi{10.1007/s10208-019-09426-y}.

\bibitem[Faradonbeh et~al.(2018)Faradonbeh, Tewari, and Michailidis]{faradonbeh2018finite}
Mohamad Kazem~Shirani Faradonbeh, Ambuj Tewari, and George Michailidis.
\newblock Finite time identification in unstable linear systems.
\newblock \emph{Automatica}, 96:\penalty0 342--353, 2018.

\bibitem[Feng and Lavaei(2021)]{feng2021learning}
Han Feng and Javad Lavaei.
\newblock Learning of dynamical systems under adversarial attacks.
\newblock In \emph{IEEE Conference on Decision and Control (CDC)}, pages 3010--3017. IEEE, 2021.

\bibitem[Feng et~al.(2023)Feng, Yalcin, and Lavaei]{feng2023learning}
Han Feng, Baturalp Yalcin, and Javad Lavaei.
\newblock Learning of dynamical systems under adversarial attacks-null space property perspective.
\newblock In \emph{American Control Conference (ACC)}, pages 4179--4184. IEEE, 2023.

\bibitem[Foster et~al.(2020)Foster, Sarkar, and Rakhlin]{foster2020learning}
Dylan Foster, Tuhin Sarkar, and Alexander Rakhlin.
\newblock Learning nonlinear dynamical systems from a single trajectory.
\newblock In \emph{Learning for Dynamics and Control}, pages 851--861. PMLR, 2020.

\bibitem[Keesman(2011)]{keesman2011system}
Karel~J Keesman.
\newblock \emph{System identification: an introduction}.
\newblock Springer Science \& Business Media, 2011.

\bibitem[Kim(2025)]{kim2025revisit}
Jihun Kim.
\newblock Revisiting the geometrically decaying step size: Linear convergence for smooth or non-smooth functions.
\newblock \emph{arXiv preprint arXiv:2508.13569}, 2025.

\bibitem[Kim and Lavaei(2025{\natexlab{a}})]{kim2024prevailing}
Jihun Kim and Javad Lavaei.
\newblock Prevailing against adversarial noncentral disturbances: Exact recovery of linear systems with the $l_1$-norm estimator.
\newblock In \emph{American Control Conference (ACC)}, pages 1161--1168. IEEE, 2025{\natexlab{a}}.

\bibitem[Kim and Lavaei(2025{\natexlab{b}})]{kim2025bridge}
Jihun Kim and Javad Lavaei.
\newblock Bridging batch and streaming estimations to system identification under adversarial attacks.
\newblock \emph{arXiv preprint arXiv:2509.15794}, 2025{\natexlab{b}}.

\bibitem[Kim and Lavaei(2025{\natexlab{c}})]{kim2025system}
Jihun Kim and Javad Lavaei.
\newblock System identification from partial observations under adversarial attacks.
\newblock \emph{arXiv preprint arXiv:2504.00244}, 2025{\natexlab{c}}.
\newblock to appear in \textit{IEEE Conference on Decision and Control}.

\bibitem[Kim et~al.(2025)Kim, Fang, and Lavaei]{kim2025sharp}
Jihun Kim, Yuchen Fang, and Javad Lavaei.
\newblock On the sharp input-output analysis of nonlinear systems under adversarial attacks.
\newblock \emph{arXiv preprint arXiv:2505.11688}, 2025.

\bibitem[Lennart(1999)]{lennart1999system}
Ljung Lennart.
\newblock System identification: theory for the user.
\newblock \emph{PTR Prentice Hall, Upper Saddle River, NJ}, 28:\penalty0 540, 1999.

\bibitem[Ljung(1978)]{ljung1978convergence}
Lennart Ljung.
\newblock Convergence analysis of parametric identification methods.
\newblock \emph{IEEE Transactions on Automatic Control}, 23\penalty0 (5):\penalty0 770--783, 1978.

\bibitem[Narendra and Annaswamy(1987)]{NarendraAnnaswamy1987}
Kumpati~S. Narendra and Anuradha~M. Annaswamy.
\newblock Persistent excitation in adaptive systems.
\newblock \emph{International Journal of Control}, 45\penalty0 (1):\penalty0 127--160, 1987.

\bibitem[Nesterov(2009)]{nesterov2009primal}
Yurii Nesterov.
\newblock Primal-dual subgradient methods for convex problems.
\newblock \emph{Mathematical Programming}, 120\penalty0 (1):\penalty0 221--259, 2009.

\bibitem[Nesterov(2013)]{nesterov2013introductory}
Yurii Nesterov.
\newblock \emph{Introductory lectures on convex optimization: A basic course}, volume~87.
\newblock Springer Science \& Business Media, 2013.

\bibitem[Polyak(1987)]{polyak1987introduction}
Boris~T. Polyak.
\newblock \emph{Introduction to Optimization}.
\newblock Optimization Software, Inc., New York, 1987.

\bibitem[Sattar and Oymak(2022)]{sattar2022non}
Yahya Sattar and Samet Oymak.
\newblock Non-asymptotic and accurate learning of nonlinear dynamical systems.
\newblock \emph{Journal of Machine Learning Research}, 23\penalty0 (140):\penalty0 1--49, 2022.

\bibitem[Schoukens and Ljung(2019)]{schoukens2019nonlinear}
Johan Schoukens and Lennart Ljung.
\newblock Nonlinear system identification: A user-oriented road map.
\newblock \emph{IEEE Control Systems Magazine}, 39\penalty0 (6):\penalty0 28--99, 2019.

\bibitem[Shor(2012)]{shor2012minimization}
Naum~Zuselevich Shor.
\newblock \emph{Minimization methods for non-differentiable functions}, volume~3.
\newblock Springer Science \& Business Media, 2012.

\bibitem[Simchowitz et~al.(2018)Simchowitz, Mania, Tu, Jordan, and Recht]{simchowitz2018learning}
Max Simchowitz, Horia Mania, Stephen Tu, Michael~I Jordan, and Benjamin Recht.
\newblock Learning without mixing: Towards a sharp analysis of linear system identification.
\newblock In \emph{Conference On Learning Theory}, pages 439--473. PMLR, 2018.

\bibitem[Tsiamis and Pappas(2019)]{tsiamis2019finite}
Anastasios Tsiamis and George~J Pappas.
\newblock Finite sample analysis of stochastic system identification.
\newblock In \emph{IEEE Conference on Decision and Control (CDC)}, pages 3648--3654. IEEE, 2019.

\bibitem[Tsiamis et~al.(2023)Tsiamis, Ziemann, Matni, and Pappas]{tsiamis2023}
Anastasios Tsiamis, Ingvar Ziemann, Nikolai Matni, and George~J. Pappas.
\newblock Statistical learning theory for control: A finite-sample perspective.
\newblock \emph{IEEE Control Systems Magazine}, 43\penalty0 (6):\penalty0 67--97, 2023.

\bibitem[Vershynin(2018)]{vershynin_2018}
Roman Vershynin.
\newblock \emph{High-Dimensional Probability: An Introduction with Applications in Data Science}.
\newblock Cambridge Series in Statistical and Probabilistic Mathematics. Cambridge University Press, 2018.

\bibitem[Wainwright(2019)]{wainwright_2019}
Martin~J. Wainwright.
\newblock \emph{High-Dimensional Statistics: A Non-Asymptotic Viewpoint}.
\newblock Cambridge Series in Statistical and Probabilistic Mathematics. Cambridge University Press, 2019.

\bibitem[Wright(1997)]{wright1997primal}
Stephen~J Wright.
\newblock \emph{Primal-dual interior-point methods}.
\newblock SIAM, 1997.

\bibitem[Yalcin et~al.(2024)Yalcin, Zhang, Lavaei, and Arcak]{yalcin2024exact}
Baturalp Yalcin, Haixiang Zhang, Javad Lavaei, and Murat Arcak.
\newblock Exact recovery for system identification with more corrupt data than clean data.
\newblock \emph{IEEE Open Journal of Control Systems}, 2024.

\bibitem[Zhang et~al.(2025)Zhang, Yalcin, Lavaei, and Sontag]{zhang2024exactrecoveryguaranteesparameterized}
Haixiang Zhang, Baturalp Yalcin, Javad Lavaei, and Eduardo~D. Sontag.
\newblock Exact recovery guarantees for parameterized nonlinear system identification problem under sparse disturbances or semi-oblivious attacks.
\newblock \emph{Transactions on Machine Learning Research}, 2025.

\bibitem[Ziemann et~al.(2022)Ziemann, Sandberg, and Matni]{ziemann2022single}
Ingvar~M Ziemann, Henrik Sandberg, and Nikolai Matni.
\newblock Single trajectory nonparametric learning of nonlinear dynamics.
\newblock In \emph{Conference on Learning Theory}, pages 3333--3364. PMLR, 2022.

\end{thebibliography}

\newpage
\section*{Appendix}

\setcounter{section}{0}
\renewcommand{\thesection}{\Alph{section}}
\renewcommand{\theHsection}{\Alph{section}} 

\section{Exact Recovery under Disturbance-Restriction}\label{app:signrestriction}

We first state the assumption on the disturbances presented in \cite{zhang2024exactrecoveryguaranteesparameterized}. 
\begin{assumption}[Semi-oblivious sub-Gaussian disturbances]\label{asp:semioblivious}
Conditional on the filtration $\mathcal{F}_t$ and the event that $\bar{d}_t\neq 0$, 
the attack vector $\bar{d}_t$ is defined by the  product $\ell_t  \hat{d}_t$, where 
\begin{enumerate}
    \item  $\hat d_t$ is a zero-mean unit vector whenever $\bar{d}_t \neq 0$, for any $\mathcal{F}_t$;
    \item $\ell_t$ is a zero-mean and sub-Gaussian with parameter $\sigma$;
    \item $\ell_t$ and $\hat d_t$ are independent whenever $\bar{d}_t \neq 0$, for any $\mathcal{F}_t$. 
\end{enumerate}
\end{assumption}
Consider $T^{(burn)}$ in Definition \ref{def:burn-in}.  \cite{zhang2024exactrecoveryguaranteesparameterized} showed that under Assumptions \ref{def:probabilistic}, \ref{persistent}, \ref{asp:stable} and \ref{asp:semioblivious}, if $T\geq T^{(burn)}$, 
then $\bar A$ is the unique solution to the estimator $\min_A \sum_{t=0}^{T-1} \|x_{t+1}-Ax_t \|_2$ with probability at least $1-\delta$ (consider $L=1$). In this section, we will show that the estimator recovers $\bar A$ for \textit{all} sufficiently large $T$ with probability at least $1-\delta$.
\begin{lemma}\label{unionall}
     Suppose that Assumptions \ref{def:probabilistic}, \ref{persistent}, \ref{asp:stable} and \ref{asp:subgaussian} hold. Consider $T^{(burn)}$ in Definition \ref{def:burn-in} and $\delta\in (0,1)$. Then, with probability at least $1-\delta$, $\bar A$ is the unique solution to the estimator $\min_A \sum_{t=0}^{T-1} \|x_{t+1}-Ax_t \|_2$ for all $T\geq T^{(\textit{burn})}$. 
\end{lemma}

\begin{proofnew}
    We first note that the claim in equations (57)-(58) of \cite{zhang2024exactrecoveryguaranteesparameterized} is also valid under Assumptions \ref{def:probabilistic}, \ref{persistent}, \ref{asp:stable} and \ref{asp:subgaussian}, with Assumption \ref{asp:semioblivious} replaced by \ref{asp:subgaussian}. \cite{zhang2024exactrecoveryguaranteesparameterized} relied on the independence between $\ell_t$ and $\hat d_t$ to show that the term of interest, $\sum_j \hat d_{k_j}^T Z g_{k_j}^{k_\ell}$ in equation (50), has zero mean. However, this independence is unnecessary, as the same result follows by applying the tower rule with respect to the filtration $\mathcal{F}_{k_j}$, using the fact that $\mathbb{E}[\hat d_{k_j} \mid \mathcal{F}_{k_j}] = 0$ under the zero-mean unit-vector assumption on $\hat d_t$ given $\mathcal{F}_t$. Thus, the sub-Gaussian parameter of the term $\sum_j \hat d_{k_j}^T Z g_{k_j}^{k_\ell}$ remains $\frac{\sigma}{1-\rho}$  (with $L=1$), which yields the same conclusion in equations (57)-(58). 

For simplicity, define $\eta:= \Theta\Bigr(\max\Bigr\{ \frac{\kappa^{10}}{(1-p)^2(1-\rho)^3 }, \frac{\kappa^4}{p(1-p)}\Bigr\}\Bigr)$. 
Equations (57)-(58) imply that the desirable event  in equation (58) at time $T$ holds with probability at least $1-\exp(-\frac{T}{\eta})$. Then, due to the union bound, the probability that such events hold for all $T\geq T^*$ is bounded below by
\begin{align}\label{union}
    1-\sum_{t=T^*}^\infty \exp\left(-\frac{t}{\eta}\right) = 1-\frac{e^{-T^*/\eta}}{1-e^{-1/\eta}} = 1-\Theta\left(\eta e^{-T^*/\eta}\right),
\end{align}
since $1-e^{-1/\eta} = \Theta(1/\eta)$. By setting $T^*:= \Theta\bigr(\eta \log(\frac{\eta}{\delta})\bigr)$,  the probability \eqref{union} is lower bounded by $1-\delta$. We substitute $\eta$ to present the required sample complexity $T^*$:
\begin{align*}
    \Theta\left( \max\left\{ \frac{\kappa^{10}}{(1-p)^2(1-\rho)^3 }, \frac{\kappa^4}{p(1-p)}\right\}\log\left(\frac{\kappa}{\delta\cdot p(1-p)(1-\rho)}\right)\right).
\end{align*}

Now, we apply the covering number argument in \cite{zhang2024exactrecoveryguaranteesparameterized}. With the $\epsilon$-covering number of the unit sphere $N = (1+\frac{2}{\epsilon})^{n^2}$ (see \cite{vershynin_2018, wainwright_2019}), we require the sample complexity of
\begin{align*}
   & \Theta\left( \max\left\{ \frac{\kappa^{10}}{(1-p)^2(1-\rho)^3 }, \frac{\kappa^4}{p(1-p)}\right\}\log\left(\frac{2N \cdot \kappa}{\delta\cdot p(1-p)(1-\rho)}\right)\right) \\ &= \Theta\left( \max\left\{ \frac{\kappa^{10}}{(1-p)^2(1-\rho)^3 }, \frac{\kappa^4}{p(1-p)}\right\}\cdot \left[\log (2N) + \log \left(\frac{  \kappa}{ p(1-p)(1-\rho)}\right) + \log\left(\frac{1}{\delta}\right)\right]\right)\\ &= \Theta\left( \max\left\{ \frac{\kappa^{10}}{(1-p)^2(1-\rho)^3 }, \frac{\kappa^4}{p(1-p)}\right\}\cdot \left[n^2 \log\left(1+\frac{2}{\epsilon}\right) + \log \left(\frac{  \kappa}{ p(1-p)(1-\rho)}\right) + \log\left(\frac{1}{\delta}\right)\right]\right)\\ &= \Theta\left( \max\left\{ \frac{\kappa^{10}}{(1-p)^2(1-\rho)^3 }, \frac{\kappa^4}{p(1-p)}\right\}\cdot \left[n^2  \log \left(\frac{  \kappa}{ p(1-p)(1-\rho)}\right) + \log\left(\frac{1}{\delta}\right)\right]\right),
\end{align*}
where the last equality comes from selecting $\log(\frac{1}{\epsilon})$ to be $\Theta\bigr(\log\bigr(\frac{\kappa}{p(1-p)(1-\rho)}\bigr)\bigr)$ (see (62)-(64) in \cite{zhang2024exactrecoveryguaranteesparameterized}). This completes the proof.
\end{proofnew}

\section{Proof of Theorem \ref{thm:angle-cosine}} \label{app:proof}

\begin{proofnew}
    Using the subgradient update iteration at time $T$, $\hat A^{(T+1)} = \hat A^{(T)} - \beta^{(T)}G_{\hat A^{(T)}, T} $, we have
    \begin{align}
        \| \hat A^{(T+1)} - \bar A \|_F^2 & = \| (\hat A^{(T)}-  \bar A) - \beta^{(T)}G_{\hat A^{(T)}, T}  \|_F^2 \notag \\
        & = \| \hat A^{(T)}-  \bar A \|_F^2 + (\beta^{(T)})^2\| G_{\hat A^{(T)}, T} \|_F^2 - 2 \beta^{(T)} \langle  G_{\hat A^{(T)}, T}, \hat A^{(T)}-  \bar A\rangle. \label{expansion}
    \end{align}
    Note that the step size $\beta^{(T)}=\frac{ \langle G_{\hat A^{(T)}, T}, \hat A^{(T)} - \bar A \rangle}{\|G_{\hat A^{(T)}, T}\|_F^2}$ proposed in the theorem minimizes \eqref{expansion}.  
    Substituting the step size yields
    \begin{align*}
        \| \hat A^{(T+1)} - \bar A \|_F^2 & = \| \hat A^{(T)}-  \bar A \|_F^2 - \frac{ (\langle G_{\hat A^{(T)}, T}, \hat A^{(T)} - \bar A \rangle)^2}{\|G_{\hat A^{(T)}, T}\|_F^2} 
         \le \sqrt{1 - \cos^2({\hat \theta_T})} \| \hat A^{(T)}- \bar A \|_F.
    \end{align*}
    The last inequality uses the  fact that $\langle G_{\hat A^{(T)}, T}, \hat A^{(T)} - \bar A \rangle = \| G_{\hat A^{(T)}, T}\|_F \| \hat A^{(T)} - \bar A\|_F \cos({\hat \theta_T})$.
\end{proofnew}

\section{Proof of Theorem \ref{thm:cosine-lower-bound}}\label{pf:cosine}

\begin{proofnew}
  Since $f_T(\cdot)$ is a sum of $\ell_2$ norm functions and is therefore convex, one can leverage subgradient inequality \citep{boyd2004convex} to derive 
  \begin{align}\label{subg1}
       f_T(\bar A)-f_T(\hat A^{(T)})  \geq \langle  \bar A - \hat{A}^{(T)},  G_{\hat{A}^{(T)}, T}  \rangle.
  \end{align}
    Thus, it holds that
    \begin{align}\label{coslowbound}
        \cos({\hat \theta_T}) = \frac{ \langle  \hat{A}^{(T)} - \bar A, G_{\hat{A}^{(T)}, T} \rangle}{\| \hat{A}^{(T)} - \bar A\|_F \| G_{\hat{A}^{(T)}, T} \|_F} 
         \ge \frac{f_T(\hat A^{(T)}) - f_T(\bar A)}{\|  \hat{A}^{(T)}-\bar A \|_F \| G_{\hat{A}^{(T)}, T} \|_F}.
    \end{align}
    The subgradient inequality also implies  
\begin{align}\label{subg2}
       f_T(\hat A^{(T)})-f_T(\bar A)  \geq \langle   \hat{A}^{(T)} - \bar A ,  G_{\bar A, T}  \rangle.
  \end{align}
  Thus, using Lemma \ref{subd_of_f_T}, we can lower bound the difference between the objectives as
    \begin{align}\label{lowboundft}
      \nonumber   f_T(\hat A^{(T)}) - f_T(\bar A)  &\ge \sup_{ G_{\bar A,T} \in \partial f_T(\bar A) } \left\{ \langle  \hat{A}^{(T)} - \bar A, G_{\bar A, T} \rangle \right\} \\
     \nonumber   & =  \sup_{ \|e_t\|\le 1, t \in \mathcal{K}^c } \left\{ \left \langle  \hat{A}^{(T)} - \bar A, -\sum_{t \in \mathcal{K}}  \hat d_t \cdot x_t^T - \sum_{t  \in \mathcal{K}^c} e_t  \cdot x_t^T \right\rangle \right\}  \\
        & =  \sum_{t  \in \mathcal{K}^c} \| ( \hat{A}^{(T)}-\bar A ) x_t \|_2-\sum_{t \in \mathcal{K}} \left\langle  (\hat{A}^{(T)}-\bar A) x_t , \frac{\bar d_t}{\| \bar d_t\|_2}\right\rangle ,
    \end{align}
    where the term on the right-hand side being positive is exactly the exact recovery condition in Corollary 3 in \cite{zhang2024exactrecoveryguaranteesparameterized} with $f(x_t) = x_t$ and $Z =  \hat{A}^{(T)}-\bar A$. Thus, from (64) in \cite{zhang2024exactrecoveryguaranteesparameterized}, along with Lemma \ref{unionall}, we obtain the relevant lemma below. 
    
    \begin{lemma} \label{lem: diff-lower-bound-haixiang}
        Suppose that  Assumptions 
   \ref{def:probabilistic}, \ref{persistent}, \ref{asp:stable} and \ref{asp:subgaussian} hold. Consider any given $B \in \mathbb{R}^{n \times n}$ and $\delta \in (0,1)$. Then,  with probability at least $1-\delta$,  the event
        \begin{align*}
         &\sum_{t  \in \mathcal{K}^c} \| (B - \bar A) x_t \|_2 -\sum_{t \in \mathcal{K}} \left\langle  B - \bar A , \frac{\bar d_t}{\| \bar d_t\|_2}\right\rangle  =\Omega \left(\frac{\lambda   p(1-p) \|B - \bar A\|_F T }{\kappa^4}\right)
    \end{align*}
    holds for all $T \ge T^{(burn)}$.
    \end{lemma}
    Using Lemma \ref{lem: diff-lower-bound-haixiang} with $B = \hat{A}^{(T)}$ and applying it to \eqref{lowboundft}, we have that, with probability at least $1-\delta$, 
    \begin{align} \label{eq:diff-lower}
         f_T(\hat A^{(T)}) - f_T(\bar A) & =\Omega\left(\frac{  \lambda p(1-p) \|\hat{A}^{(T)} - \bar A\|_F T }{\kappa^4}\right)
    \end{align}
   holds for all $T\geq T^{(burn)}$. In other words, for a sufficiently large $T$, $f_T(\cdot)$ has a sharpness with respect to $\bar A$; \textit{i.e.}, $\exists \bar c>0$ such that $f_T(\hat A^{(T)})-f_T(\bar A) \geq \bar c \|\hat A^{(T)}- \bar A\|_F$. 
   
   Next, we upper-bound the norm of the subgradient as below:
    \begin{align}
          \| G_{\hat{A}^{(T)}, T} \|_F^2 
         & = \Big\langle - \sum_{t=0}^{T-1} \partial \norm{(\bar A -\hat A^{(T)}) x_t + \bar d_t}_2  \cdot x_t^T,  - \sum_{t=0}^{T-1} \partial \norm{(\bar A -\hat A^{(T)}) x_t + \bar d_t}_2  \cdot x_t^T \Big \rangle \notag  \\
         & \le \sum_{t = 0 }^{T-1} \sum_{m = 0}^{T-1} \|x_t \|_2 \|x_m \|_2 = \left( \sum_{t = 0}^{T-1} \|x_t\|_2 \right)^2. \label{eq:gradient-upper}
    \end{align}
    The inequality holds due to the norm of subgradients being at most $1$ and the Cauchy-Schwarz inequality. 
    In the next lemma, we upper bound $\sum_{t=0}^{T-1}\|x_t\|_2$ with high probability. 
\begin{lemma}\label{sumxt}
Suppose that Assumptions \ref{asp:stable} and \ref{asp:subgaussian} hold and consider $\delta \in (0,1)$. Then, with probability at least $1-\delta$, we have 
\begin{align}\label{sumxtexp}
    \sum_{t=0}^{T-1} \|x_t\|_2 =\mathcal{O}\left(\frac{\sigma T}{1-\rho}\right)
\end{align}
for all $T=\Omega\bigr(\log(\frac{1}{\delta})\bigr)$.
\end{lemma}
We defer the proof of Lemma \ref{sumxt} to the end of this section. 
     We now consider \eqref{eq:diff-lower} and \eqref{sumxtexp} to hold with probability at least $1-\frac{\delta}{2}$. Since this does not affect the order of the sample complexity, applying the union bound implies that both events hold for all $T\geq T^{(burn)}$ with probability at least $1-\delta$ (note that we already have $T^{(burn)}=\Omega(\log(\frac{1}{\delta})) $). 
     
     We substitute \eqref{eq:diff-lower}, \eqref{eq:gradient-upper}, \eqref{sumxtexp} into \eqref{coslowbound}  to arrive at a lower bound on $\cos({\hat \theta_T})$:
\begin{align}\label{forpolyak}
     \nonumber \cos({\hat \theta_T})
        & \ge \frac{f_T(\hat A^{(T)}) - f_T(\bar A)}{\|  \hat{A}^{(T)}-\bar A \|_F \| G_{\hat{A}^{(T)}, T} \|_F}  =\Omega\left(\frac{\lambda  p(1-p) \|\hat A^{(T)}-\bar A\|_F  T}{\kappa^4  } \cdot \frac{1}{ \|\hat A^{(T)}-\bar A\|_F \cdot \frac{\sigma}{1-\rho}T} \right)
         \\&\hspace{37mm} =\Omega\left(\frac{\lambda  p(1-p) (1-\rho)  }{\sigma \kappa^4  }  \right)=\Omega\left(\frac{  p(1-p) (1-\rho)  }{ \kappa^5  }  \right)>0
\end{align}
for all $T\geq T^{(burn)}$ with probability at least $1-\delta$.  This completes the proof.
\end{proofnew}
\textbf{Proof of Lemma \ref{sumxt}~}
\textit{Due to system dynamics \eqref{sysdyn}, we have 
    \begin{align}\label{submul}
       \nonumber \sum_{t=0}^{T-1} \|x_t\|_2 &= \sum_{t=0}^{T-1} \Bigr\|\bar A^t x_0 + \sum_{i=0}^{t-1} \bar A^{t-1-i} \bar d_{i}\Bigr\|_2 < \sum_{i=0}^{\infty} \rho^i   \Bigr[\|x_0\|_{2} + \sum_{t=0}^{T-2} \|\bar d_t\|_{2}\Bigr]\\&=\frac{1}{1-\rho}\Bigr[\|x_0\|_{2} + \sum_{t=0}^{T-2} \|\bar d_t\|_{2}   \Bigr] =\frac{1}{1-\rho}\Bigr[\|x_0\|_{2} + \sum_{t=0}^{T-2} \ell_t   \Bigr] 
    \end{align}
    due to the triangle inequality.
    Under Assumption \ref{asp:subgaussian}, the sub-Gaussian parameter of $\ell_t$ is $\sigma$.
The centering lemma in \cite{vershynin_2018} implies that $\ell_t - \mathbb{E}[\ell_t]$ is a zero-mean sub-Gaussian variable with parameter $\mathcal{O}(\sigma)$. It follows that 
\begin{align*}
    \mathbb{E}\left[\exp\left(\nu \sum_{t=0}^{T-2} (\ell_t - \mathbb{E}[\ell_t])\right) \right]&= \mathbb{E}\left[\mathbb{E}\left[\exp\left(\nu \sum_{t=0}^{T-2} (\ell_t - \mathbb{E}[\ell_t])\right) ~\Biggr|~\mathcal{F}_{T-2} \right]\right]\\&= \mathbb{E}\left[\exp\left(\nu \sum_{t=0}^{T-3} (\ell_t - \mathbb{E}[\ell_t])\right) \mathbb{E}\Bigr[\exp\left(\nu  (\ell_{T-2} - \mathbb{E}[\ell_{T-2}])\right) ~|~\mathcal{F}_{T-2} \Bigr]\right]\\&\le \mathbb{E}\left[\exp\left(\nu \sum_{t=0}^{T-3} (\ell_t - \mathbb{E}[\ell_t])\right) \right] \exp(\nu^2 \mathcal{O}(\sigma^2)) \\&\le \cdots \le  \exp(\nu^2 (T-1) \mathcal{O}(\sigma^2)) 
\end{align*}
for all $\nu\in\mathbb{R}$, which implies that $\sum_{t=0}^{T-2} \ell_t - \mathbb{E}[\ell_t]$ is a zero-mean sub-Gaussian variable with parameter $\mathcal{O}(\sqrt{T}\sigma)$, and thus $\|x_0\|_2 - \mathbb{E}[\|x_0\|_2] + \sum_{t=0}^{T-2} \ell_t - \mathbb{E}[\ell_t]$ also has a sub-Gaussian parameter of $\mathcal{O}(\sqrt{T} \sigma)$.}

\textit{We can leverage Hoeffding's inequality to arrive at   
    \begin{align*}
        \mathbb{P}\Bigr( \|x_0\|_2 - \mathbb{E}[\|x_0\|_2] + \sum_{t=0}^{T-2} \ell_t - \mathbb{E}[\ell_t]\leq \sigma T  \Bigr) \geq 1-\exp\Bigr(-\Omega\Bigr(\frac{\sigma^2 T^2}{\sigma^2 T}\Bigr)\Bigr) = 1-\exp(-\Omega(T)).
    \end{align*}
    Note that the sub-Gaussian parameter of $\ell_t$ being $\sigma$ yields  $\mathbb{E}[\|x_0\|_2] + \sum_{t=0}^{T-2} \mathbb{E}[\ell_t] =\mathcal{O}(\sigma T)$. This concludes that \begin{align}\label{hoeff}
        \mathbb{P}\Bigr( \|x_0\|_2 + \sum_{t=0}^{T-2} \ell_t \leq \mathcal{O}(\sigma T)  \Bigr) \geq 1-\exp(-\Omega(T)).
    \end{align} 
Due to the union bound, the probability that  \eqref{hoeff} holds for all $T\geq \tilde T$ is bounded below by
\begin{align}\label{temp}
   1- \sum_{t=\tilde T}^\infty \exp(-cT)  =   1-\frac{e^{-c\tilde T}}{1-e^{-c}} 
\end{align}
for some $c>0$. Setting $\tilde T := \frac{1}{c}\log(\frac{1}{\delta(1-e^{-c})}) = \Theta\bigr(\log (\frac{1}{\delta})\bigr)$ makes \eqref{temp} equal to $1-\delta$, which implies that  the event in \eqref{hoeff} holds for all $T\geq \tilde T$ with probability at least $1-\delta$. Finally, considering   \eqref{submul} completes the proof.  }

\section{Proof of Theorem \ref{thm:polyak-result}}\label{pf:polyak}

\begin{proofnew}
  We first note that $f_T(\hat A^{(T)})-f_T(\bar A)\ge 0 $ when $T\ge T^{(burn)}$ with high probability due to Lemma \ref{unionall}. It follows that the proposed step size in the theorem satisfies $\beta^{(T)}\ge 0$. 
    
    Now, using the subgradient update iteration at time $T$, $\hat A^{(T+1)} = \hat A^{(T)} - \beta^{(T)}G_{\hat A^{(T)}, T} $, we have
    \begin{align}
        \| \hat A^{(T+1)} - \bar A \|_F^2 & = \| (\hat A^{(T)}-  \bar A) - \beta^{(T)}G_{\hat A^{(T)}, T}  \|_F^2 \notag \\
        & = \| (\hat A^{(T)}-  \bar A) \|_F^2 + (\beta^{(T)})^2\| G_{\hat A^{(T)}, T} \|_F^2 - 2 \beta^{(T)} \langle  G_{\hat A^{(T)}, T}, \hat A^{(T)} -  \bar A\rangle \notag \\
        & \le \| (\hat A^{(T)}-  \bar A) \|_F^2 + (\beta^{(T)})^2\| G_{\hat A^{(T)}, T} \|_F^2  - 2 \beta^{(T)} (f_T(\hat A^{(T)}) - f_T(\bar A)), \label{expansion_new}
    \end{align}
    where the last inequality leverages \eqref{subg1} and $\beta^{(T)}\ge 0$. Substituting $\beta^{(T)}$ yields
    \begin{align*}
        \| \hat A^{(T+1)} - \bar A \|_F^2 & \le  \| \hat A^{(T)}-  \bar A \|_F^2 - \frac{ (f_T(\hat A^{(T)}) - f_T(\bar A) )^2}{\|G_{\hat A^{(T)}, T}\|_F^2} \\& =  \| \hat A^{(T)}-  \bar A \|_F^2 \Biggr(1 - \Biggr(\frac{ f_T(\hat A^{(T)}) - f_T(\bar A) }{\| \hat A^{(T)}-  \bar A \|_F\|G_{\hat A^{(T)}, T}\|_F}\Biggr)^2\Biggr),
    \end{align*}
    where the right-hand side can be bounded above by $\gamma^2$ with high probability  when $T\ge T^{(burn)}$, which follows from \eqref{forpolyak}.
     %
     Thus, we have convergence to the ground-truth system with the linear rate $\gamma$ since 
     \begin{align*}
         \| \hat A^{(T+1)} - \bar A \|_F & \le  \gamma\| \hat A^{(T)}-  \bar A\|_F. 
     \end{align*}
     The rest of the arguments are the same to those in Corollary \ref{cor:linear-convergence-best}.
\end{proofnew}

\section{Proof of Theorems \ref{thm:constant-step-size} and \ref{thm: diminishing-step-size}}\label{app:constant-diminishing}
\begin{proofnew} We first define the following lemma and defer its proof to the end of this section.

    \begin{lemma} \label{lem:objective-difference}
    Let $D = \| \hat A^{(T^{(burn)})} - \bar A \|_F$ be the distance of the iterate to the ground-truth $\bar A$ at the end of the burn-in time $T^{(burn)}$. Then, for $\bar T \ge T^{(burn)}$, it holds that
    \begin{align}\label{generalmin}
      \min_{T\in [T^{(burn)},\bar T]} \|\hat{A}^{(T)} - \bar A\|_F =\mathcal{O}\left(\frac{ \kappa^4 \bigr[D^2 + (\frac{\sigma }{1-\rho} )^2\sum_{T=T^{(burn)}}^{\bar T} (\beta^{(T)}T)^2 \bigr]}{ \lambda p(1-p)   \sum_{T=T^{(burn)}}^{\bar T}  (\beta^{(T)} T)}\right)
    \end{align}
    with probability at least $1-\delta$.
\end{lemma}

 For Theorem \ref{thm:constant-step-size}, setting the step size to be $\beta^{(T)} = \beta$ for all $T\ge T^{(burn)}$ implies the term on the right-hand side to be
    \begin{align}\label{constantmin}
       \mathcal{O}\left(\frac{ \kappa^4 \bigr[D^2 + (\frac{\sigma }{1-\rho} )^2 \beta^2 \sum_{t=T^{(burn)}}^{\bar T} t^2\bigr] }{ \lambda p(1-p) \beta   \sum_{t=T^{(burn)}}^{\bar T}  t}\right).
    \end{align}
    Minimizing the right-hand side with respect to $\beta$ results in the following constant step size: 
    \begin{align*}
        \beta^{(T)} = \beta =  \Theta\left(\frac{D(1-\rho)}{\sigma  (\sum_{t=T^{(burn)}}^{\bar T} t^2)^{1/2}}\right).
    \end{align*}
    Substituting this constant step size back into \eqref{constantmin}  provides the simplified upper bound that holds with probability at least $1 - \delta$:
    \begin{align*}
        \min_{T\in [T^{(burn)},\bar T]} \|\hat{A}^{(T)} - \bar A\|_F  &=\mathcal{O}\left(\frac{\kappa^4}{\lambda p(1-p)} \cdot \frac{D\sigma}{1-\rho} \cdot \frac{(\sum_{t=T^{(burn)}}^{\bar T} t^2)^{1/2}}{\sum_{t=T^{(burn)}}^{\bar T} t}\right)\\&=\mathcal{O}\left(\frac{D\kappa^5}{p(1-p)(1-\rho)}  \cdot \frac{(\sum_{t=T^{(burn)}}^{\bar T} t^2)^{1/2}}{\sum_{t=T^{(burn)}}^{\bar T} t}\right).
    \end{align*}
    This completes the proof of Theorem \ref{thm:constant-step-size}. Now, to prove Theorem \ref{thm: diminishing-step-size}, 
    setting the step size to be $\beta^{(T)} = \beta/T$ implies the right-hand side of  \eqref{generalmin} to be
\begin{align}\label{dimmin}
    \mathcal{O}\left(\frac{ \kappa^4 \bigr[D^2 + (\frac{\sigma }{1-\rho} )^2\beta^2 (\bar T-T^{(burn)})\bigr]}{ \lambda p(1-p)     \beta (\bar T-T^{(burn)})}\right).
\end{align}
Minimizing the right-hand side with respect to $\beta$ results in the following constant:
\begin{align*}
        \beta := \Theta \left( \frac{D  (1-\rho)}{\sigma  (\bar T- T^{(burn)})^{1/2}} \right).
    \end{align*}
    Substituting this constant back into \eqref{dimmin} ensures the following upper bound with probability at least $1-\delta$:
    \begin{align*}
         \min_{T\in [T^{(burn)},\bar T]} \|\hat{A}^{(T)} - \bar A\|_F &=\mathcal{O}\left(\frac{\kappa^4}{\lambda p(1-p)}\cdot \frac{D\sigma}{1-\rho} \cdot \frac{(\bar T-T^{(burn)})^{1/2}}{\bar T-T^{(burn)}}\right)\\ &= \mathcal{O}\left(\frac{D\kappa^5}{p(1-p)(1-\rho)}  \cdot \frac{1}{(\bar T-T^{(burn)})^{1/2}}\right).
    \end{align*}
This completes the proof of Theorem \ref{thm: diminishing-step-size}.
\end{proofnew}
 \textbf{Proof of Lemma \ref{lem:objective-difference}~}  
 \textit{We note that the proof hinges on the sharpness condition given in \eqref{eq:diff-lower}. 
 Using \eqref{expansion_new} and the iterative substitution between the time periods from $T^{(burn)}$ to $\bar T$, we obtain
    \begin{align*}
         \| &\hat A^{(\bar T+1)} - \bar A \|_F^2 \\ & \le 
         \| \hat A^{(T^{(burn)})} -  \bar A\|_F^2+ \sum_{T=T^{(burn)}}^{\bar T} (\beta^{(T)})^2\| G_{\hat A^{(T)}, T} \|_F^2 
        - 2\sum_{T=T^{(burn)}}^{\bar T}  \beta^{(T)} (f_T(\hat A^{(T)}) - f_T(\bar A)) \\
        & \le  D^2 + \sum_{T=T^{(burn)}}^{\bar T} (\beta^{(T)} )^2 \left(\frac{\sigma T }{1-\rho} \right)^2  - 2\sum_{T=T^{(burn)}}^{\bar T}  (\beta^{(T)}) \cdot \Omega\left(\frac{  \lambda p(1-p) \|\hat{A}^{(T)} - \bar A\|_F  \cdot T}{\kappa^4}\right).
    \end{align*}
    The last inequality follows from  the fact that $\| G_{\hat A^{(T)}, T} \|_F$ is bounded above by $\frac{\sigma T}{1-\rho}$ considering \eqref{eq:gradient-upper}-\eqref{sumxtexp},  and $f_T(\hat A^{(T)})-f_T(\bar A)$ bounded below by \eqref{eq:diff-lower} with high probability for all $T\ge T^{(burn)}$.  Since $ \| \hat A^{(\bar T+1)} - \bar A \|_F^2 $ is nonnegative, we arrive at 
\begin{align*}
    2\sum_{T=T^{(burn)}}^{\bar T}  (\beta^{(T)} T)\cdot \Omega\left(\frac{  \lambda p(1-p) \|\hat{A}^{(T)} - \bar A\|_F }{\kappa^4}\right)\leq  D^2 + \sum_{T=T^{(burn)}}^{\bar T} (\beta^{(T)}T)^2 \left(\frac{\sigma }{1-\rho} \right)^2 ,
\end{align*}
where the left-hand side can further be bounded below by $$2 \Omega\left(\frac{  \lambda p(1-p)  }{\kappa^4}\right) \cdot \min_{T\in [T^{(burn)},\bar T]} \|\hat{A}^{(T)} - \bar A\|_F \cdot \sum_{T=T^{(burn)}}^{\bar T}  (\beta^{(T)} T).$$ Thus, it follows that
\begin{align*}
    \min_{T\in [T^{(burn)},\bar T]} \|\hat{A}^{(T)} - \bar A\|_F  =\mathcal{O}\left(\frac{ D^2 + (\frac{\sigma }{1-\rho} )^2\sum_{T=T^{(burn)}}^{\bar T} (\beta^{(T)}T)^2 }{\frac{  \lambda p(1-p)  }{\kappa^4}  \sum_{T=T^{(burn)}}^{\bar T}  (\beta^{(T)} T)}\right)
\end{align*}
holds with high probability. This completes the proof.}

\end{document}